 \documentclass[12pt]{article}

 \usepackage{mathrsfs,amsfonts}
 \usepackage[dvips]{color}

 \newtheorem{theorem}{Theorem}[section]
 \newtheorem{lemma}[theorem]{Lemma}
 \newtheorem{corollary}[theorem]{Corollary}
 \newtheorem{proposition}[theorem]{Proposition}
 \newtheorem{example}[theorem]{Example}
 \newtheorem{remark}[theorem]{Remark}

 \def\blemma{\begin{lemma}\sl{}\def\elemma{\end{lemma}}}
 \def\bproposition{\begin{proposition}\sl{}\def\eproposition{\end{proposition}}}
 \def\btheorem{\begin{theorem}\sl{}\def\etheorem{\end{theorem}}}
 \def\bcorollary{\begin{corollary}\sl{}\def\ecorollary{\end{corollary}}}
 \def\bremark{\begin{remark}\sl{}\def\eremark{\end{remark}}}

 \def\mcr{\mathscr}\def\mbb{\mathbb}\def\mbf{\mathbf}

 \def\sign{\mathrm{sign}}

 \def\<{\langle}\def\>{\rangle}\def\ar{\!\!&}\def\nnm{\nonumber}

 \def\proof{\noindent{\it Proof.~~}}\def\qed{\hfill$\Box$\medskip}

\begin{document}


\bigskip\bigskip

\centerline{\Large\bf Strong solutions for stochastic differential}

\medskip

\centerline{\Large\bf equations with jumps}

\bigskip

\centerline{Zenghu Li\footnote{Supported by NSFC (10525103 and
10721091) and CJSP.}}

\centerline{School of Mathematical Sciences, Beijing Normal
University,}

\centerline{Beijing 100875, People's Republic of China}

\centerline{E-mail: \tt lizh@bnu.edu.cn}

\medskip

\centerline{Leonid Mytnik\footnote{Supported by Israel Science Foundation
grant No.1162/06.}}

\centerline{Faculty of Industrial Engineering and Management,}

\centerline{Technion---Israel Institute of Technology, Haifa 32000,
Israel}

\centerline{E-mail: \tt leonid@ie.technion.ac.il}

\bigskip

{\narrower{\narrower

\noindent\textit{Abstract.} General stochastic equations with jumps are
studied. We provide criteria for the uniqueness and existence of strong
solutions under non-Lipschitz conditions of Yamada-Watanabe type. The
results are applied to stochastic equations driven by spectrally positive
L\'evy processes.

\smallskip

\noindent\textit{Mathematics subject classifications (2000).}
Primary: 60H10, 60H20; Secondary: 60J80.

\smallskip

\noindent\textit{Keywords:} Stochastic equation, strong solution,
pathwise uniqueness, non-Lipschitz condition.

\par}\par}

\bigskip


\section{Introduction}

\setcounter{equation}{0}

The question of pathwise uniqueness for one-dimensional stochastic
differential equations driven by one-dimensional Brownian motions has
been resolved a long time ago by Yamada and Watanabe~\cite{YW71}; see
also Barlow~\cite{B82}. The same question can also be asked for
stochastic differential equations driven by discontinuous L\'evy noises.
Let us consider the equation
 \begin{eqnarray}\label{1.1}
dx(t)= F(x(t-))dL_t, \qquad t\geq 0.
 \end{eqnarray}
Bass~\cite{B03} and Komatsu~\cite{K82} showed that if $\{L_t\}$ is a
symmetric stable process with exponent $\alpha\in (1,2)$ and if $x\mapsto
F(x)$ is a bounded function with modulus of continuity $z\mapsto \rho(z)$
satisfying
 \begin{eqnarray}\label{1.2}
\int_{0+} \frac{1}{\rho(z)^\alpha}dz = \infty,
 \end{eqnarray}
then (\ref{1.1}) admits a strong solution and the solution is pathwise
unique. This condition is the analogue of the Yamada-Watanabe criterion
for the diffusion coefficient. In particular, if $F$ is H\"older
continuous with exponent $1/\alpha$, then the pathwise uniqueness holds
for (\ref{1.1}). The required H\"older exponent tends to $1/2$ as
$\alpha\to 2$ and it tends to $1$ (Lipschitz condition) as $\alpha\to 1$.
When the integral in (\ref{1.2}) is finite, Bass~\cite{B03} constructed a
continuous function $x\mapsto \phi(x)$ having continuity modulus
$x\mapsto \rho(x)$ for which the pathwise uniqueness for (\ref{1.1})
fails; see also~\cite{BBC04}.

The pathwise uniqueness and strong solutions for stochastic differential
equations driven by spectrally positive L\'evy noises were studied in
\cite{FL10}. Those equations arise naturally in the study of branching
processes. A typical special continuous state branching process is the
non-negative solution to the stochastic differential equation
 \begin{eqnarray}\label{1.3}
dx(t) = \sqrt[\alpha]{x(t-)} dL_t, \qquad t\ge 0,
 \end{eqnarray}
where $\{L_t\}$ is a Brownian motion (for $\alpha= 2$) or a spectrally
positive $\alpha$-table process (for $1<\alpha< 2$). Note that the
coefficient $x\mapsto \sqrt[\alpha]{x}$ in (\ref{1.3}) is non-decreasing,
non-Lipschitz and degenerate at the origin. More general stochastic
equations with similar structures arise naturally in limit theorems of
branching processes with interactions or/and immigration.

In this paper we consider a class of stochastic differential equations
with jumps, which generalizes the equation (\ref{1.3}). This exploration
can be regarded as a continuation of~\cite{FL10}. We extend the results
of~\cite{FL10} in two directions. First of all, we notice that the
pathwise uniqueness results proved in~\cite{FL10} for non-negative
c\`adl\`ag solutions can easily be extended to any c\`adl\`ag solutions.
This extended result is given in Proposition~\ref{pr3.1}. Its proof,
which is in fact the most involved stochastic part behind the results in
this paper, goes through along the same lines as in~\cite{FL10}.

The second direction is to apply the above result to formulate some
criteria for the pathwise uniqueness and existence of strong solutions to
general stochastic differential equations with jumps. We consider this to
be the main part of this paper. The proofs in this part involve some
analytical arguments that allow us to apply the general pathwise
uniqueness criterion of Proposition~\ref{pr3.1}. From those results we
derive sufficient conditions for the existence and uniqueness of
non-negative strong solutions under suitable additional assumptions.

We also give applications of our main results to stochastic equations
driven by spectrally positive L\'evy processes. These extend and improve
substantially the results of~\cite{FL10}. As a consequence of one of
those results we get the following counterpart of the theorem of Bass
\cite{B03}:

 \btheorem\label{th1.1}
Let $\{L_t\}$ be a spectrally positive stable process with exponent
$\alpha\in (1,2)$ that is, there exists $c_{\alpha}$ such that
 \begin{eqnarray*}
\mbf{E}\left[e^{-uL(t)}\right]
 =
e^{-c_{\alpha}u^{\alpha}t}, \qquad t\geq 0, u\geq 0.
 \end{eqnarray*}
Let $F$ be a non-decreasing function on $\mbb{R}$ with modulus of
continuity $z\mapsto \rho(z)$ satisfying
 \begin{eqnarray}\label{1.4}
\int_{0+} \frac{1}{\rho(z)^{\alpha/(\alpha-1)}}\,dz =\infty.
 \end{eqnarray}
Also assume that there is a constant $K\ge 0$ such that
 \begin{eqnarray*}
|F(x)| \leq K(1+|x|),
\quad x\in \mbb{R}.
 \end{eqnarray*}
Then there is a pathwise unique strong solution to (\ref{1.1}). \etheorem

By the above theorem, if $F$ is a non-decreasing function H\"older
continuous with exponent $1-1/\alpha$, then the pathwise uniqueness holds
for (\ref{1.1}). The required H\"older exponent tends to $0$ as
$\alpha\to 1$, which differs sharply from the criterion of
Bass~\cite{B03} for a symmetric stable noise. Note that this result is
also consistent with the Yamada-Watanabe result in the sense that as
$\alpha\rightarrow 2$ the critical H\"older exponent converges to $1/2$.

The organization of the paper is as follows. The main theorem is stated
in Section~\ref{sec:2}. Its proof is provided in Section~\ref{sec:3}. In
Section~\ref{sec:4} a number of particular cases is considered, for
example, SDE's with stable L\'evy noises. Theorem~\ref{th1.1} is a
consequence of one of the results obtained in that section. Throughout
this paper, we make the conventions
 \begin{eqnarray*}
\int_a^b = \int_{(a,b]}
 \quad\mbox{and}\quad
\int_a^\infty = \int_{(a,\infty)} \qquad\mbox{for}\quad  b\ge a\in \mbb{R}.
 \end{eqnarray*}


\section{Main strong uniqueness and existence results}
\label{sec:2}

\setcounter{equation}{0}

Suppose that $\mu_0(du)$ and $\mu_1(du)$ are $\sigma$-finite measures on
the complete separable metric spaces $U_0$ and $U_1$, respectively. Let
$(\Omega, \mcr{G}, \mcr{G}_t, \mbf{P})$ be a filtered probability space
satisfying the usual hypotheses. Let $\{B(t)\}$ be a standard
$(\mcr{G}_t)$-Brownian motion and let $\{p_0(t)\}$ and $\{p_1(t)\}$ be
$(\mcr{G}_t)$-Poisson point processes on $U_0$ and $U_1$ with
characteristic measures $\mu_0(du)$ and $\mu_1(du)$, respectively.
Suppose that $\{B(t)\}$, $\{p_0(t)\}$ and $\{p_1(t)\}$ are independent of
each other. Let $N_0(ds,du)$ and $N_1(ds,du)$ be the Poisson random
measures associated with $\{p_0(t)\}$ and $\{p_1(t)\}$, respectively.
Suppose in addition that
 \begin{itemize}

\item $x\mapsto \sigma(x)$ is a continuous function on $\mbb{R}$;

\item $x\mapsto b(x)$ is a continuous function on $\mbb{R}$ having
    the decomposition $b=b_1-b_2$ with $b_2$ being continuous and
    non-decreasing;

\item $(x,u) \mapsto g_0(x,u)$ is a  Borel function on
    $\mbb{R} \times U_0$ such that $x \mapsto
    g_0(x,u)$ is non-decreasing for every $u\in U_0$;

\item $(x,u) \mapsto g_1(x,u)$ is a Borel function on $\mbb{R}
    \times U_1$.

 \end{itemize}
Let $\tilde{N}_0(ds,du)$ be the compensated measure of $N_0(ds,du)$. By a
\textit{solution of} the stochastic equation
 \begin{eqnarray}\label{2.1}
x(t)
 \ar=\ar
x(0) + \int_0^t \sigma(x(s-))dB(s) + \int_0^t\int_{U_0} g_0(x(s-),u)
\tilde{N}_0(ds,du) \nnm \\
 \ar \ar
+ \int_0^t b(x(s-))ds + \int_0^t\int_{U_1} g_1(x(s-),u) N_1(ds,du)
 \end{eqnarray}
we mean a c\`adl\`ag and $({\mcr{G}}_t)$-adapted real process $\{x(t)\}$
that satisfies the equation almost surely for every $t\ge0$. Since $x(s-)
\neq x(s)$ for at most countably many $s\ge0$, we can also use $x(s)$
instead of $x(s-)$ for the integrals with respect to $dB(s)$ and $ds$ on
the right hand side of (\ref{2.1}). We say \textit{pathwise uniqueness}
holds for (\ref{2.1}) if for any two solutions $\{x_1(t)\}$ and
$\{x_2(t)\}$ of the equation satisfying $x_1(0) = x_2(0)$ we have $x_1(t)
= x_2(t)$ almost surely for every $t\ge 0$. Let $(\mcr{F}_t)_{t\ge 0}$ be
the augmented natural filtration generated by $\{B(t)\}$, $\{p_0(t)\}$
and $\{p_1(t)\}$. A solution $\{x(t)\}$ of (\ref{2.1}) is called a
\textit{strong solution} if $x(t)$ is measurable with respect to
$\mcr{F}_t$ for every $t\ge 0$; see~\cite[p.163]{IW89} or
\cite[p.76]{S05}.

\blemma\label{le2.1} Suppose that $(z\land z^2)\nu(dz)$ is a finite
measure on $(0,\infty)$ and define
 \begin{eqnarray}\label{2.2}
\alpha_\nu = \inf\Big\{\beta>1: ~ \lim_{x\to 0+} \, x^{\beta-1} \!
\int_x^\infty z\nu(dz) = 0\Big\}.
 \end{eqnarray}
Then $1\le \alpha_\nu\le 2$ and, for any $\alpha> \alpha_\nu$,
 \begin{eqnarray}\label{2.3}
\lim_{x\to 0+} \, x^{\alpha-2} \! \int_0^x z^2 \nu(dz) = 0.
 \end{eqnarray}
\elemma

\proof By (\ref{2.2}) it is clear that $\alpha_\nu\ge 1$. For $x>0$
let
 \begin{eqnarray*}
G(x) = \int_x^\infty z \nu(dz)
 \quad\mbox{and}\quad
H(x) = \int_0^x z^2 \nu(dz).
 \end{eqnarray*}
Given $\varepsilon>0$, choose $a>0$ so that $H(a)<\varepsilon$. Then
for $a\ge x>0$ we have
 \begin{eqnarray*}
xG(x) = x\int_x^a z \nu(dz) + xG(a) \le \int_x^a z^2 \nu(dz) + xG(a) \le
\varepsilon + xG(a).
 \end{eqnarray*}
It follows that $\limsup_{x\to 0+} xG(x) \le \varepsilon$. That proves
$\lim_{x\to 0+} xG(x) = 0$, and so $\alpha_\nu\le 2$. Clearly,
(\ref{2.3}) holds for any $\alpha\ge 2$. By integration by parts,
 \begin{eqnarray}\label{2.4}
H(x) = - \int_0^x z dG(z) = - xG(x) + \int_0^x G(z)dz.
 \end{eqnarray}
Thus we have
 \begin{eqnarray*}
\lim_{x\to 0+} \int_0^x G(z)dz = \lim_{x\to 0+} H(x) + \lim_{x\to
0+}xG(x) = 0.
 \end{eqnarray*}
Now suppose that $\alpha_\nu< \alpha< 2$. In view of (\ref{2.2}), for any
$\varepsilon>0$ there exists $b>0$ so that $x^{\alpha-1}G(x)<
\varepsilon$ for all $0<x\le b$. Then (\ref{2.4}) implies
 \begin{eqnarray*}
x^{\alpha-2}H(x) \le x^{\alpha-2}\int_0^x G(z)dz \le x^{\alpha-2}
\int_0^x \varepsilon z^{1-\alpha} dz = \varepsilon (2-\alpha)^{-1},
 \end{eqnarray*}
and hence $\lim_{z\to 0+} x^{\alpha-2}H(x) = 0$. \qed

Let us consider a set $U_2\subset U_1$ satisfying $\mu_1(U_1\setminus
U_2)<\infty$. As in the proof of Proposition~2.2 in~\cite{FL10} one can
show that the uniqueness/existence of strong solutions for (\ref{2.1})
can be reduced to the same question for the equation with $U_1$ replaced
by $U_2$. Then in what follows all conditions for the ingredients of
(\ref{2.1}) only involve $U_2$ instead of $U_1$. As usual, let us
consider some growth conditions on the coefficients:
 \begin{itemize}

\item[{\rm(2.a)}] There is a constant $K\ge 0$ such that
    \begin{eqnarray*}
    \ar\ar \sigma(x)^2 + \int_{U_0} g_0(x,u)^2\mu_0(du)
    + \int_{U_2} g_1(x,u)^2\mu_1(du) \cr
    \ar\ar\qquad
    +\, b(x)^2 + \bigg(\int_{U_2} |g_1(x,u)| \mu_1(du)\bigg)^2\leq K(1+x^2),
    \quad x\in \mbb{R}.
    \end{eqnarray*}

 \end{itemize}

We next introduce our main conditions on the modulus of continuity that
are particularly useful in applications to stochastic equations driven by
L\'evy processes. The conditions are given as follows:
 \begin{itemize}

\item[{\rm(2.b)}] For each $m\ge 1$ there is a non-decreasing and
    concave function $z\mapsto r_m(z)$ on $\mbb{R}_+$ such that
    $\int_{0+} r_m(z)^{-1}\, dz=\infty$ and \begin{eqnarray*}
    |b_1(x)-b_1(y)| + \int_{U_2} |l_1(x,y,u)| \mu_1(du)\le r_m(|x-y|)
    \end{eqnarray*} for $|x|,|y|\le m$, where
    $l_1(x,y,u)=g_1(x,u)-g_1(y,u)$.

\item[{\rm(2.c)}] For each $m\ge 1$ there is a constant $p_m>0$, a
    non-decreasing function $z\mapsto \rho_m(z)$ on $\mbb{R}_+$ and a
    function $u\mapsto f_m(u)$ on $U_0$ such that \begin{eqnarray*}
    \int_{0+} \rho_m(z)^{-2}dz = \infty, \qquad \int_{U_0}
    [f_m(u)\land f_m(u)^2]\mu_0(du)<\infty \end{eqnarray*} and
    \begin{eqnarray*} |\sigma(x)-\sigma(y)| \le \rho_m(|x-y|), \quad
    |g_0(x,u)-g_0(y,u)|\le \rho_m(|x-y|)^{2p_m}f_m(u)
    \end{eqnarray*}\bigskip for all $|x|,|y|\le m$ and $u\in U_0$.

 \end{itemize}
For each $m\ge 1$ and the function $f_m$ defined in (2.c) we define the
constant
 \begin{eqnarray*}
\alpha_m := \inf\Big\{\beta>1: ~ \lim_{x\to 0+} \, x^{\beta-1} \!
\int_{U_0} f_m(u)1_{\{f_m(u)\ge x\}}\mu_0(du) = 0\Big\}.
 \end{eqnarray*}
By Lemma~\ref{le2.1} we have $1\le \alpha_m\le 2$. Our first main theorem
of this paper is the following

\btheorem\label{th2.2} Suppose that conditions (2.a,b,c) hold with
 \begin{eqnarray}\label{2.5}
p_m>1 - 1/\alpha_m ~ \mbox{for} ~ \alpha_m<2,
 ~~\mbox{or}~
p_m=1/2 ~ \mbox{for} ~ \alpha_m=2.
 \end{eqnarray}
Then for any given $x(0)\in \mbb{R}$, there exists a pathwise unique
strong solution $\{x(t)\}$ to (\ref{2.1}).
 \etheorem

From the above theorem we may derive some results on non-negative
solutions of (\ref{2.1}). For that purpose let us consider the following
conditions:
 \begin{itemize}

 \item[{\rm(2.d)}] $\sigma(0) = 0$, $b(0)\ge 0$ and $g_0(0,u) = 0$
     for $u\in U_0$, and $g_1(x,u)+x\ge 0$ for $x\in\mbb{R}_+$ and
     $u\in U_1$;

 \item[{\rm(2.e)}] There is a constant $K\ge 0$ such that
 \begin{eqnarray*} b(x) + \int_{U_2} |g_1(x,u)| \mu_1(du) \le K(1+x), \qquad x\ge
0;
 \end{eqnarray*}

 \item[{\rm(2.f)}] There is a non-decreasing function $x\mapsto L(x)$
     on $\mbb{R}_+$ so that
 \begin{eqnarray*} \sigma(x)^2 + \int_{U_0} [|g_0(x,u)|\land g_0(x,u)^2]
\mu_0(du) \le L(x), \qquad x\ge 0.
 \end{eqnarray*}

 \end{itemize}
By Proposition~2.1 of~\cite{FL10}, under condition (2.d) any solution of
(\ref{2.1}) with non-negative initial value remains non-negative forever.

\btheorem\label{th2.3} Suppose that conditions (2.b,c,d,e,f) hold with
(\ref{2.5}). Then for any given $x(0)\in \mbb{R}_+$, there exists a
pathwise unique non-negative strong solution $\{x(t)\}$ to
(\ref{2.1}). \etheorem

 \bremark\label{re2.4}
Under the conditions of Theorem~\ref{th2.3} we can actually conclude that
for any given $x(0)\in \mbb{R}_+$ there is a pathwise unique strong
solution to (\ref{2.1}) and the solution is non-negative. That follows
from Proposition~2.1 of \cite{FL10}.
 \eremark

 \bremark\label{re2.5}
Note that when $\alpha_m<2$ the assumptions of Theorem~\ref{th2.2}
and~\ref{th2.3} are strictly weaker than Theorems~2.5 and~5.3
of~\cite{FL10}. In some particular cases the condition~(\ref{2.5}) can be
weakened to $p_m \geq 1 - 1/\alpha_m$, as in the case of stable driving
noise. This is done in Theorem~\ref{th4.2}.
 \eremark

\section{Proofs of Theorems~\ref{th2.2} and \ref{th2.3}}\label{sec:3}

\setcounter{equation}{0}

The crucial part of the proof of Theorem~\ref{th2.2} is verifying the
pathwise uniqueness for (\ref{2.1}). As we have mentioned already it is
enough to consider the equation
 \begin{eqnarray}\label{3.1}
x(t)
 \ar=\ar
x(0) + \int_0^t \sigma(x(s-))dB(s) + \int_0^t\int_{U_0} g_0(x(s-),u)
\tilde{N}_0(ds,du) \nnm \\
 \ar \ar
+ \int_0^t b(x(s-))ds + \int_0^t\int_{U_2} g_1(x(s-),u) N_1(ds,du)
 \end{eqnarray}
For a function $f$ defined on the real line $\mbb{R}$, note
 \begin{eqnarray*}
\Delta_zf(x) = f(x+z) - f(x)
 \quad\mbox{and}\quad
D_zf(x) = \Delta_zf(x) - f^\prime(x)z.
 \end{eqnarray*}
We shall need the next result, which provides a criterion for the
pathwise uniqueness. It extends the criterion of Theorem~3.1
in~\cite{FL10}, where it was formulated just for non-negative solutions.

 \bproposition\label{pr3.1}
Suppose that condition (2.b,c) holds. Then the pathwise uniqueness of
solution to~(\ref{3.1}) holds if for each $m\geq 1$ there exists a
sequence of non-negative and twice continuously differentiable functions
$\{ \phi_k\}$ with the following properties:
 \begin{itemize}

\item[(i)] $\phi_k(z)\mapsto |z|$ non-decreasingly as $k\rightarrow
    \infty$;

\item[(ii)] $0\leq \phi_k'(z)\leq 1$ for $z\geq 0$ and $-1\leq
    \phi_k'(z)\leq 0$ for $z\leq 0$;

\item[(iii)] $\phi_k''(z)\geq 0$ for $z\in \mbb{R}$ and as
    $k\rightarrow \infty$,
 \begin{eqnarray*}
\phi_k''(x-y) [\sigma(x)-\sigma(y)]^2\rightarrow 0
 \end{eqnarray*}
uniformly on $|x|,|y|\leq m$;

\item[(iv)] as $k\rightarrow \infty$,
 \begin{eqnarray*}
\int_{U_0} D_{l_0(x,y,u)} \phi_k(x-y)\mu_0(du) \rightarrow 0
 \end{eqnarray*}
uniformly on $|x|,|y|\leq m$, where $l_0(x,y,u)=g_0(x,u)-g_0(y,u)$.

 \end{itemize}
 \eproposition

\proof For non-negative solutions the result was given in Theorem~3.1
of~\cite{FL10}. In what follows we will show that the proof
in~\cite{FL10} goes through for any c\`adl\`ag solutions. Let
$\{x_1(t)\}$ and $\{x_2(t)\}$ be any two solutions of (\ref{3.1})
starting at $x_1(0)=x_2(0) = x_0$. For each $m\geq 1$ define
$\tau_m=\inf\{t\geq 0:\; |x_1(t)|\geq m$ or $|x_2(t)|\geq m\}$. Recall
that $l_i(x,y,u) = g_i(x,u)-g_i(y,u)$, $i=0,1$. By (\ref{3.1}) and the
It\^o formula one can show
 \begin{eqnarray*}
\phi_k(\zeta(t\wedge \tau_m))
 \ar=\ar
\int_0^{t\wedge \tau_m} \phi_k'(\zeta(s-)) [b(x_1(s-)) -
b(x_2(s-))]\,ds \cr
 \ar\ar
+\, \frac{1}{2}\int_0^{t\wedge \tau_m} \phi_k''(\zeta(s-))
[\sigma(x_1(s-))-\sigma(x_2(s-))] ds \cr
 \ar\ar
+ \int_0^{t\wedge \tau_m} \,ds \int_{U_2} \Delta_{l_1(x_1(s-),x_2(s-),u)}
\phi_k(\zeta(s-))\mu_1(du) \cr
 \ar\ar
+ \int_0^{t\wedge \tau_m} \,ds \int_{U_0} D_{l_0(x_1(s-),x_2(s-),u)}
\phi_k(\zeta(s-))\mu_0(du) \nnm \\
 \ar\ar
+\, M_m(t),
 \end{eqnarray*}
where
 \begin{eqnarray*}
M_m(t)
 \ar=\ar
\int_0^{t\land\tau_m}\phi_k^\prime(\zeta(s-))
[\sigma(x_1(s-)) - \sigma(x_2(s-))]dB(s) \cr
 \ar \ar
+ \int_0^{t\land\tau_m}\int_{U_2} \Delta_{l_1(x_1(s-),x_2(s-),u)} \phi_k(\zeta(s-))
\tilde{N}_1(ds,du) \cr
 \ar \ar
+ \int_0^{t\land\tau_m}\int_{U_0} \Delta_{l_0(x_1(s-),x_2(s-),u)}\phi_k(\zeta(s-))
\tilde{N}_0(ds,du).
 \end{eqnarray*}
Under conditions (2.b,c) it is easy to show that $\{M_m(t)\}$ is a
martingale. Therefore, we can follow the same argument as in the proof of
Theorem~3.1 of~\cite{FL10} to get that, as $k\rightarrow \infty$,
 \begin{eqnarray*}
\mbf{E}[|\zeta(t\wedge\tau_m)|]
 \leq
\int_0^{t} r_m(\mbf{E}[|\zeta(s\wedge\tau_m)|]) ds.
 \end{eqnarray*}
From this by standard argument we have $\mbf{E}[|\zeta(t\wedge \tau_m)|]
= 0$ for every $t\geq 0$. Since $\{x_1(t)\}$ and $\{x_2(t)\}$ are
c\`adl\`ag, we have that $\tau_m\rightarrow \infty$ as
$m\rightarrow\infty$. Hence letting $m\rightarrow \infty$ and using the
right continuity of $\{\zeta(t)\}$ we get the result. \qed

To prove the pathwise uniqueness for (\ref{3.1}) we need to introduce
more notation and prove a lemma which will play a crucial role in the
proofs. For each integer $m\ge1$ we shall construct a sequence of
functions $\{\phi_k\}$ that satisfies the properties required in
Proposition~\ref{pr3.1}. Although main ideas are similar to those in the
proof of Theorem~3.2 of~\cite{FL10}, we will go through the details for
the sake of completeness. Let $1=a_0> a_1> a_2> \ldots> 0$ be defined by
 \begin{eqnarray*}
\int_{a_k}^{a_{k-1}}\rho_m(z)\,dz = k.
 \end{eqnarray*}
Let $x\mapsto \psi_k(x)$ be a non-negative continuous function on
$\mbb{R}$ satisfying $\int_{a_k}^{a_{k-1}} \psi_k(x)\,dx=1$ and
 \begin{eqnarray}\label{3.2}
0\leq \psi_k(x)
 \leq
2k^{-1} \rho_m(x)^{-2}1_{(a_k, a_{k-1})}(x).
 \end{eqnarray}
For each $k\geq 1$ we define the non-negative and twice continuously
differentiable function
 \begin{eqnarray*}
\phi_k(z)= \int_0^{|z|}\,dy \int_0^{y}\psi_k(x)\,dx, \quad z\in
\mbb{R}.
 \end{eqnarray*}
Note that although the sequences $\{a_k\}$, $\{\phi_k\}$ and $\{\psi_k\}$
also depend on $m\ge 1$, we do not put this additional index to simplify
the notation.

\blemma\label{le3.2} Suppose that condition (2.c) holds. Fix $m\ge 1$ and
let $a_k$, $\phi_k$ and $\psi_k$ be defined as above. Then the sequence
$\{\phi_k\}$ satisfies properties (i)--(iii) in Proposition~\ref{pr3.1}
and for any $h> 0$,
\begin{eqnarray}\label{3.3}
\lefteqn{\int_{U_0} D_{l_0(x,y,u)}\phi_k(x-y) \mu_0(du)}
\nnm\\
 \ar\le\ar
k^{-1} \rho_m(|x-y|)^{4p_m -2}1_{\{|x-y|\le a_{k-1}\}} \int_{U_0}
f_m(u)^21_{\{f_m(u)\le h\}} \mu_0(du) \nnm \\
 \ar \ar
+\, \rho_m(|x-y|)^{2p_m}1_{\{|x-y|\le a_{k-1}\}}\int_{U_0} f_m(u)
1_{\{f_m(u)>h\}}\mu_0(du).
 \end{eqnarray}
\elemma

\proof By definition, the sequence $\{\phi_k\}$ satisfies properties (i)
and (ii) in Proposition~\ref{pr3.1}. Moreover, by (\ref{3.2}) we get
 \begin{eqnarray}\label{3.4}
\phi_k^{\prime\prime}(x)
 =
\psi_k(|x|)
 \le
2k^{-1}\rho_m(|x|)^{-2}1_{(a_k, a_{k-1})}(|x|)
 \end{eqnarray}
for all $x\in\mbb{R}$. This together with condition (2.c) implies
 \begin{eqnarray*}
\phi_k^{\prime\prime}(x-y)[\sigma(x) - \sigma(y)]^2
 \le
\psi_k(|x-y|)\rho_m(|x-y|)^2
 \le
2/k
 \end{eqnarray*}
for $|x|,|y|\le m$. Thus $\{\phi_k\}$ also satisfies property (iii) in
Proposition~\ref{pr3.1}. Observe that
 \begin{eqnarray}\label{3.5}
D_z\phi_k(x-y)
 =
\Delta_z\phi_k(x-y) - \phi_k^\prime(x-y)z
 \le
|z|1_{\{|x-y|\le a_{k-1}\}}
 \end{eqnarray}
when $(x-y)z\ge 0$. By Taylor's expansion,
 \begin{eqnarray*}
D_z\phi_k(x-y)
 =
z^2\int_0^1\phi_k^{\prime\prime}(x-y+tz)(1-t)dt
 =
z^2\int_0^1\psi_k(|x-y+tz|)(1-t)dt.
 \end{eqnarray*}
Then (\ref{3.4}) and the monotonicity of $\zeta\mapsto \rho_m(\zeta)$ imply
 \begin{eqnarray}\label{3.6}
D_z\phi_k(x-y)
 \ar\le\ar
2k^{-1}z^2\int_0^1\frac{(1-t)1_{(a_k,a_{k-1})}(|(x-y)+tz|)}
{\rho_m(|(x-y)+tz|)^2}dt \nnm \\
 \ar\le\ar
k^{-1}z^2\rho_m(|x-y|)^{-2}1_{\{|x-y|\le a_{k-1}\}}
 \end{eqnarray}
when $(x-y)z\ge 0$ and $|x|,|y|\le m$. Recall that $l_0(x,y,u) =
g_0(x,u)-g_0(y,u)$. Since $x\mapsto g_0(x,u)$ is non-decreasing, for
$|x|,|y| \le m$ we get by (\ref{3.5}) and (2.c) that
 \begin{eqnarray*}
D_{l_0(x,y,u)}\phi_k(x-y)
 \le
|l_0(x,y,u)|1_{\{|x-y|\le a_{k-1}\}}
 \le
\rho_m(|x-y|)^{2p_m} f_m(u)1_{\{|x-y|\le a_{k-1}\}}.
 \end{eqnarray*}
Similarly, by (\ref{3.6}) and (2.c) we have
 \begin{eqnarray*}
D_{l_0(x,y,u)}\phi_k(x-y)
 \ar\le\ar
k^{-1}\rho_m(|x-y|)^{-2}l_0(x,y,u)^21_{\{|x-y|\le a_{k-1}\}} \nnm \\
 \ar\le\ar
k^{-1}\rho_m(|x-y|)^{4p_m -2}f_m(u)^21_{\{|x-y|\le a_{k-1}\}}.
 \end{eqnarray*}
Then (\ref{3.3}) follows immediately. \qed

\bproposition\label{pr3.3} Under the conditions (2.b,c) and (\ref{2.5}),
the pathwise uniqueness holds for equation (\ref{3.1}). \eproposition

\proof For $\alpha_m=2$ and $p_m=2$, the result was essentially proved in
Theorem~3.3 of~\cite{FL10} for non-negative solutions. It follows along
the same lines for all solutions. So we here only consider the case of
$\alpha_m<2$ and $p_m>1 - 1/\alpha_m$. By Lemma~\ref{le3.2} we get that
the sequence $\{\phi_k\}$ satisfies properties (i)---(iii) in
Proposition~\ref{pr3.1}. Moreover for any $\beta>0$ we can take $h=
\rho_m(|x-y|)^{2\beta}$ in (\ref{3.3}) to get
 \begin{eqnarray*}
\lefteqn{\int_{U_0} D_{l_0(x,y,u)}\phi_k(x-y) \mu_0(du)} \cr
 \ar\le\ar
k^{-1} \rho_m(|x-y|)^{2(2p_m-1)}1_{\{|x-y|\le a_{k-1}\}} \int_{U_0}
f_m(u)^21_{\{f_m(u)\le\rho_m(|x-y|)^{2\beta}\}} \mu_0(du) \cr
 \ar \ar
+\, \rho_m(|x-y|)^{2p_m}1_{\{|x-y|\le a_{k-1}\}}\int_{U_0} f_m(u)
1_{\{f_m(u)>\rho_m(|x-y|)^{2\beta}\}} \mu_0(du).
 \end{eqnarray*}
Since $\lim_{k\to \infty}a_k=0$ and $\lim_{z\to 0+} \rho(z)=0$, for
$\alpha_m< \alpha< 2$ we use Lemma~\ref{le2.1} to see
 \begin{eqnarray}\label{3.7}
\lefteqn{\int_{U_0} D_{l_0(x,y,u)}\phi_k(x-y) \mu_0(du)} \cr
 \ar\le\ar
k^{-1} \rho_m(|x-y|)^{2(2p_m-1)}\rho_m(|x-y|)^{2\beta(2-\alpha)}
1_{\{|x-y|\le a_{k-1}\}} \cr
 \ar \ar
+\, \rho_m(|x-y|)^{2p_m}\rho_m(|x-y|)^{2\beta(1-\alpha)}
1_{\{|x-y|\le a_{k-1}\}}
 \end{eqnarray}
when $k\ge 1$ is sufficiently large. If we can choose $\beta$ and
$\alpha$ in the way that
 \begin{eqnarray*}
2(2p_m-1) + 2\beta(2-\alpha)> 0
 \quad\mbox{and}\quad
2p_m + 2\beta(1-\alpha)> 0,
 \end{eqnarray*}
the value on the right hand side of (\ref{3.7}) will tend to zero as $k\to
\infty$. The requirement is equivalent to
 \begin{eqnarray*}
\frac{1-2p_m}{2-\alpha} < \beta < \frac{p_m}{\alpha-1},
 \end{eqnarray*}
which can be done as long as
 \begin{eqnarray*}
\frac{1-2p_m}{2-\alpha} < \frac{p_m}{\alpha-1}
 \end{eqnarray*}
or, equivalently, $p_m>1-1/\alpha$. For that purpose it sufficient to have
$p_m> 1 - 1/\alpha_m$. This gives property (iv) in Proposition~\ref{pr3.1}
and hence the pathwise uniqueness for (\ref{3.1}). \qed

\bproposition\label{pr3.4} Suppose that conditions (2.a) hold. Let
$\{x(t)\}$ be a solution of (\ref{3.1}) with $\mbf{E}[x(0)^2]<\infty$.
Then we have
 \begin{eqnarray}\label{3.8}
\mbf{E}\Big[1+\sup_{0\le s\le t}x(s)^2\Big]
 \le
(1+6\mbf{E}[x(0)^2])\exp\{6K(4+t)t\}.
 \end{eqnarray}
\eproposition

\proof Let $\tau_m = \inf\{t\ge0: |x(t)|\ge m\}$ for $m\ge 1$. Since
$\{x(t)\}$ has c\`adl\`ag sample paths, we have $\tau_m\to \infty$ as
$m\to \infty$. Let us rewrite (\ref{3.1}) into
 \begin{eqnarray*}
x(t)
 \ar=\ar
x(0) + \int_0^t \sigma(x(s-))dB(s) + \int_0^t\int_{U_0} g_0(x(s-),u)
\tilde{N}_0(ds,du) \nnm \\
 \ar \ar
+ \int_0^t b(x(s-))ds + \int_0^t\int_{U_2} g_1(x(s-),u)
\tilde{N}_1(ds,du) \nnm \\
 \ar \ar
+ \int_0^tds\int_{U_2} g_1(x(s-),u) \mu_1(du).
 \end{eqnarray*}
By Doob's martingale inequalities we have
 \begin{eqnarray*}
\mbf{E}\Big[\sup_{0\le s\le t}x(s\land\tau_m)^2\Big]
 \ar\le\ar
6\mbf{E}[x(0)^2] + 24\mbf{E}\bigg[\int_0^{t\land\tau_m}\sigma(x(s-))^2ds\bigg] \\
 \ar \ar
+\, 6\mbf{E}\bigg[\bigg(\int_0^{t\land\tau_m} |b(x(s-))| ds\bigg)^2\bigg] \\
 \ar \ar
+\, 24\mbf{E}\bigg[\int_0^{t\land\tau_m}ds\int_{U_0}g_0(x(s-),u)^2\mu_0(du)\bigg] \\
 \ar \ar
+\, 24\mbf{E}\bigg[\int_0^{t\land\tau_m}ds\int_{U_2}g_1(x(s-),u)^2\mu_1(du)\bigg] \\
 \ar \ar
+\, 6\mbf{E}\bigg[\bigg(\int_0^{t\land\tau_m}ds\int_{U_2} |g_1(x(s-),u)|
\mu_1(du)\bigg)^2\bigg] \\
 \ar\le\ar
6\mbf{E}[x(0)^2] + 6K(4+t)\mbf{E}\bigg[\int_0^{t\land\tau_m} (1+x(s-)^2) ds\bigg].
 \end{eqnarray*}
Then it is easy to see that
 $$
t\mapsto F_m(t) := \mbf{E}\Big[\sup_{0\le s\le t}x(s\land\tau_m)^2\Big]
 $$
is locally bounded on $[0,\infty)$. Since $s\mapsto x(s)$ has at most a
countable number of jumps, from the above inequality we obtain
 \begin{eqnarray*}
1+F_m(t)
 \ar\le\ar
1+6\mbf{E}[x(0)^2] + 6K(4+t)\mbf{E}\bigg[\int_0^{t\land\tau_m} (1+x(s)^2) ds\bigg] \cr
 \ar\le\ar
1+6\mbf{E}[x(0)^2] + 6K(4+t)\int_0^t [1+F_m(s)] ds.
 \end{eqnarray*}
By Gronwall's inequality,
 \begin{eqnarray*}
\mbf{E}\Big[1+\sup_{0\le s\le t}x(s\land\tau_m)^2\Big]
 \ar\le\ar
(1+6\mbf{E}[x(0)^2])\exp\{6K(4+t)t\}.
 \end{eqnarray*}
Then (\ref{3.8}) follows by Fatou's lemma. \qed

\noindent{\textit{Proof of Theorem~\ref{th2.2}~} \textit{Step~1)} Suppose
that conditions (2.b,c) and (\ref{2.5}) hold. Instead of condition (2.a),
we here assume there is a constant $K\ge 0$ such that
 \begin{eqnarray}\label{3.9}
\ar\ar \sigma(x)^2 + b(x)^2 + \sup_{u\in U_0}|g_0(x,u)| + \int_{U_0}
g_0(x,u)^2 \mu_0(du)
 \cr
 \ar\ar\qquad
+ \int_{U_2}g_1(x,u)^2 \mu_1(du) + \bigg(\int_{U_2} |g_1(x,u)|
\mu_1(du)\bigg)^2\leq K,
\quad x\in \mbb{R}.
 \end{eqnarray}
Let $\{V_n\}$ be a non-decreasing sequence of Borel subsets of $U_0$ so
that $\cup_{n=1}^\infty V_n = U_0$ and $\mu_0(V_n)<\infty$ for every
$n\ge 1$. By the result on continuous-type stochastic equations, there is
a weak solution to
 \begin{eqnarray}\label{3.10}
x(t)
 \ar=\ar
x(0) + \int_0^t \sigma(x(s))dB(s) + \int_0^t
b(x(s))ds \nnm \\
 \ar \ar
- \int_0^tds\int_{V_n} g_0(x(s),u)\mu_0(du);
 \end{eqnarray}
see, e.g., Ikeda and Watanabe (1989, p.169). By Proposition~\ref{pr3.3},
the pathwise uniqueness holds for (\ref{3.10}), thus the equation has a
pathwise unique strong solution. Let $\{W_n\}$ be a non-decreasing
sequence of Borel subsets of $U_1$ so that $\cup_{n=1}^\infty W_n = U_2$
and $\mu_1(W_n)<\infty$ for every $n\ge 1$. Then for every integer $n\ge
1$ there is a unique strong solution $\{x_n(t): t\ge0\}$ to
 \begin{eqnarray*}
x(t)
 \ar=\ar
x(0) + \int_0^t \sigma(x(s-))dB(s) + \int_0^t\int_{V_n} g_0(x(s-),u)
\tilde{N}_0(ds,du) \nnm \\
 \ar \ar
+ \int_0^t b(x(s-))ds + \int_0^t\int_{W_n} g_1(x(s-),u) N_1(ds,du).
 \end{eqnarray*}
As in the proof of Lemma~4.3 of~\cite{FL10} one can see the sequence
$\{x_n(t)\}$ is tight in $D([0,\infty),\mbb{R})$. Following the proof of
Theorem~4.4 of~\cite{FL10} it is easy to show that any limit point of the
sequence is a weak solution to (\ref{3.1}). This and
Proposition~\ref{pr3.3} imply the existence and uniqueness of the strong
solution to (\ref{3.1}); see, e.g.,~\cite[p.104]{S05}.

\textit{Step~2)} Suppose that the original conditions (2.a,b,c) and
(\ref{2.5}) hold. For each $m\ge 1$ let
 \begin{eqnarray*}
\chi_m(x)&=&\left\{\begin{array}{lcl}
                  x,&{\rm if}& |x|\leq m,\\
                  m, &{\rm if}& x> m,\\
                  -m, &{\rm if}& x< -m.\\
                 \end{array}\right.
 \end{eqnarray*}
We consider the equation
 \begin{eqnarray}\label{3.11}
x(t)
 \ar=\ar
x(0) + \int_0^t \sigma(\chi_m(x(s-)))dB(s) + \int_0^t
b_m(\chi_m(x(s-)))ds \nnm \\
 \ar \ar
+ \int_0^t\int_{U_0} \chi_m\circ g_0(\chi_m(x(s-)),u)
\tilde{N}_0(ds,du) \nnm \\
 \ar \ar
+ \int_0^t\int_{U_2} g_1(\chi_m(x(s-)),u) N_1(ds,du),
 \end{eqnarray}
where
 \begin{eqnarray*}
b_m(x) = b(x) - \int_{U_0} [g_0(x,u) - \chi_m\circ g_0(x,u)] \mu_0(du).
 \end{eqnarray*}
By the first step, there is a unique strong solution to (\ref{3.11}).
Then using Proposition~\ref{pr3.4} one can show as in the proof of
Proposition~2.4 of~\cite{FL10} that there is a pathwise unique strong
solution to (\ref{3.1}). Hence as we have mentioned above, there is a
pathwise unique strong solution to (\ref{2.1}) (see Proposition~2.2
of~\cite{FL10} and its proof for analogous result). \qed

\noindent{\textit{Proof of Theorem~\ref{th2.3}~} By Proposition~2.1 of
\cite{FL10} and Theorem~\ref{th2.2} there is a pathwise unique
non-negative strong solution $\{x_m(t)\}$ to the equation
 \begin{eqnarray}
\label{3.12}
x(t)
 \ar=\ar
x(0) + \int_0^t \sigma\left(\chi_m(x(s-)\vee 0)\right)dB(s) + \int_0^t
b(\chi_m(x(s-)\vee 0)) ds \nnm  \\
 \ar \ar
+ \int_0^t\int_{U_0} \chi_m\circ g_0(\chi_m(x(s-)\vee 0),u)
\tilde{N}_0(ds,du) \nnm \\
 \ar \ar
+ \int_0^t\int_{U_2} \chi_m\circ g_1(x(s-)\vee 0,u) N_1(ds,du).
 \end{eqnarray}
By Proposition~2.3 of~\cite{FL10} the first moment of $\{x_m(t)\}$ is
dominated by a locally bounded function on $[0,\infty)$ independent of
$m\ge 1$. Then one can follow the proof of Proposition~2.4 of~\cite{FL10}
to show there is a pathwise unique non-negative strong solution to
 \begin{eqnarray}\label{3.13}
x(t)
 \ar=\ar
x(0) + \int_0^t \sigma(x(s-)\vee 0)dB(s) + \int_0^t
b(x(s-)\vee 0) ds \nnm \\
\nonumber
 \ar \ar
+ \int_0^t\int_{U_0}  g_0((x(s-)\vee 0),u)
\tilde{N}_0(ds,du) \nnm \\
 \ar \ar
+ \int_0^t\int_{U_2} g_1(x(s-)\vee 0,u) N_1(ds,du).
 \end{eqnarray}
Now note that the non-negative solution to~(\ref{3.13}) is also the
non-negative  solution  to~(\ref{3.1}). This and Proposition~\ref{pr3.3}
imply that there is a pathwise unique non-negative strong solution to
(\ref{3.1}). This again as in the proof of Theorem~\ref{th2.2} implies
that there is a pathwise unique non-negative strong solution to
(\ref{2.1}). \qed

 \bremark\label{re3.5}
The above proofs show it is unnecessary to assume the existence of the
sequence $\{V_n\}$ in (4.b) of~\cite{FL10}. As a consequence, condition
(5.b) of~\cite{FL10} is also unnecessary for the results in Section~5 of
that paper.
 \eremark


\section{Stochastic equations with L\'evy noises}\label{sec:4}

\setcounter{equation}{0}

In this section, we give some applications of our main results to
stochastic equations driven by L\'evy processes. Let $(\sigma,b)$ be
given as in Section~2 and let $\nu_0(dz)$ and $\nu_1(dz)$ be
$\sigma$-finite Borel measures on $(0,\infty)$ satisfying
 \begin{eqnarray*}
\int_0^\infty (z\land z^2) \nu_0(dz) + \int_0^\infty (1\land z)
\nu_1(dz)< \infty.
 \end{eqnarray*}
Let $\alpha_0$ be the constant defined by (\ref{2.2}) for the
measure $\nu_0(dz)$. In addition, we suppose that
 \begin{itemize}

\item $x\mapsto h_0(x)$ is a continuous and non-decreasing function
    on $\mbb{R}$;

\item $x\mapsto h_1(x)$ is a continuous function on $\mbb{R}$.

 \end{itemize}
Suppose we have a filtered probability space $(\Omega, {\mcr{G}},
{\mcr{G}}_t, \mbf{P})$ satisfying the usual hypotheses. Let $\{B(t)\}$ be
an $({\mcr{G}}_t)$-Brownian motion and let $\{L_0(t)\}$ and $\{L_1(t)\}$
be $({\mcr{G}}_t)$-L\'evy processes with exponents
 \begin{eqnarray*}
u\mapsto \int_0^\infty \big(e^{iuz} - 1 - iuz\big) \nu_0(dz)
 \quad\mbox{and}\quad
u\mapsto \int_0^\infty \big(e^{iuz} - 1\big) \nu_1(dz),
 \end{eqnarray*}
respectively. Suppose that $\{B(t)\}$, $\{L_0(t)\}$ and $\{L_1(t)\}$ are
independent of each other. Note that $\{L_0(t)\}$ is centered and
$\{L_1(t)\}$ is non-decreasing. We introduce the conditions:
 \begin{itemize}

\item[{\rm(4.a)}] There is a constant $K\ge 0$ such that
 \begin{eqnarray*}
|\sigma(x)| + |b(x)| + |h_0(x)| + |h_1(x)|\leq K(1+|x|),
\quad x\in \mbb{R};
 \end{eqnarray*}

\item[{\rm(4.b)}] There exists a non-decreasing and concave function
    $z \mapsto r(z) $ on $\mbb{R}_+$ such that
    $\int_{0+}r(z)^{-1}\,dz =\infty$ and
 \begin{eqnarray*}
|b(x)-b(y)| + |h_1(x)-h_1(y)| \leq r(|x-y|), \qquad x,y\in \mbb{R};
 \end{eqnarray*}

\item[{\rm(4.c)}] There is a constant $p>0$ and a non-decreasing
    function $z\mapsto \rho(z)$ on $\mbb{R}_+$ such that
    $\int_{0+}\rho (z)^{-2}\,dz = \infty$ and
 \begin{eqnarray*}
|\sigma(x)-\sigma(y)| + |h_0(x)-h_0(y)|^{1/2p} \leq \rho(|x-y|),
\qquad x,y \in \mbb{R};
 \end{eqnarray*}

\item[{\rm(4.d)}] $\sigma(0) = h_0(0) = 0$, $b(0)\ge 0$, and
    $h_1(x)\ge 0$ for $x\in\mbb{R}_+$;

\item[{\rm(4.e)}] There is a constant $K\ge 0$ such that
 \begin{eqnarray*}
b(x)+h_1(x)\leq K(1+x), \qquad x\ge 0.
 \end{eqnarray*}

  \end{itemize}

\btheorem\label{th4.1} {\rm(i)} If conditions (4.a,b,c) are satisfied
with $p> 1 - 1/\alpha_0$, then for any given $x(0)\in \mbb{R}$ there is a pathwise unique strong solution to
 \begin{eqnarray}\label{4.1}
dx(t) = \sigma(x(t))dB(t) + h_0(x(t-))dL_0(t) + b(x(t))dt +
h_1(x(t-))dL_1(t).
 \end{eqnarray}
{\rm(ii)} If conditions (4.b,c,d,e) are satisfied with $p> 1 -
1/\alpha_0$, then for any given $x(0)\in \mbb{R}_+$ there is a pathwise
unique non-negative strong solution to (\ref{4.1}). \etheorem

\proof
By L\'evy-It\^o decompositions, the L\'evy processes have the
following representations
 \begin{eqnarray*}
L_0(t) \ar=\ar \int_0^t\int_0^1 z \tilde{N}_0(ds,dz) - \int_0^tds\int_1^\infty z \nu_0(dz) \\
 \ar\ar
+ \int_0^t\int_1^\infty z N_1(ds,dz,\{0\}), \\
L_1(t) \ar=\ar \int_0^t\int_0^\infty z N_1(ds,dz,\{1\}),
 \end{eqnarray*}
where $N_0(ds,dz)$ and $N_1(ds,dz,du)$ are Poisson random
measures with intensities
 \begin{eqnarray*}
1_{\{z\leq 1\}}ds\nu_0(dz) \quad\mbox{and}\quad ds[1_{\{z>1\}}\nu_0(dz)\delta_0(du)+\nu_1(dz)\delta_1(du)],
 \end{eqnarray*}
respectively, and $\tilde{N}_0(ds,dz)$ is the compensated measure of
$N_0(ds,dz)$. Here $N_0(ds,dz)$ and $N_1(ds,dz,du)$ are independent and they are independent of $\{B(t)\}$. By applying Theorem~\ref{th2.2} with
 \begin{eqnarray*}
U_0=(0,1], ~ U_1 = [(1,\infty)\times \{0\}]\cup[(0,\infty)\times \{1\}] \quad\mbox{and}\quad U_2=(0,1]\times \{1\},
 \end{eqnarray*}
we see that there is a pathwise unique strong
solution to
 \begin{eqnarray*}
x(t) \ar=\ar x(0) + \int_0^t\sigma(x(s)) dB(s) +
\int_0^t\int_0^1 h_0(x(s-)) z \tilde{N}_0(ds,dz) \nnm \\
 \ar \ar
+ \int_0^t \left(b(x(s))-h_0(x(s))\int_1^{\infty}z\nu_0(dz)\right)ds \nnm \\
 \ar \ar
+ \int_0^t\int_{U_1} g_1(x(s-),z,u) N_1(ds,dz,du),
 \end{eqnarray*}
where
 \begin{eqnarray*}
g_1(x,z,u) = h_0(x)z1_{\{z>1,u=0\}} + h_1(x)z1_{\{u=1\}}.
 \end{eqnarray*}
However, this is just another form of the equation (\ref{4.1}) and hence
part (i) of the theorem follows. The proof of part (ii) is similar.  \qed

 \btheorem\label{th4.2}
Suppose that $\{B(t)\}$, $\{L_0(t)\}$ and $\{L_1(t)\}$ are given as the
above with $\nu_0(dz) = z^{-1-\alpha}dz$ for $1<\alpha<2$.  Then we have:\\
{\rm(i)} If conditions (4.a,b,c) are satisfied with $p\ge 1 - 1/\alpha$,
then for any given $x(0)\in \mbb{R}$ there is a pathwise unique strong
solution to (\ref{4.1});\\ {\rm(ii)} If conditions (4.b,c,d,e) are
satisfied with $p\ge 1 - 1/\alpha$, then for any given $x(0)\in
\mbb{R}_+$ there is a pathwise unique non-negative strong solution to
(\ref{4.1}).
 \etheorem

\proof Let $\{a_k\}$, $\{\phi_k\}$ and $\{\psi_k\}$ be defined as before
Lemma~\ref{le3.2} with $\rho_m=\rho$. Then we can easily apply
Lemma~\ref{le3.2} to get that $\{\phi_k\}$ satisfies properties (i)-(iii)
in Proposition~\ref{pr3.1}. Moreover, using again Lemma~\ref{le3.2} with
$\mu_0(du)=u^{-1-\alpha}du$, $p_m=p = (\alpha-1)/\alpha$, $\rho_m=\rho$
and $f_m(u)=u$ we can rewrite~(\ref{3.3}) as
 \begin{eqnarray*}
\lefteqn{\int_{0}^{\infty} D_{l_0(x,y,u)}\phi_k(x-y) \mu_0(du)} \\
 \ar\le\ar
k^{-1} \rho(|x-y|)^{4p -2} \int_{0}^h u^{1-\alpha}\,du
+ \rho(|x-y|)^{2p}\int_h^{\infty} u^{-\alpha}\, du  \\
 \ar=\ar
k^{-1}(2-\alpha)^{-1} \rho(|x-y|)^{4(\alpha-1)/\alpha-2}
h^{2-\alpha} + (\alpha-1)^{-1} \rho(|x-y|)^{2(\alpha-1)/\alpha} h^{1-\alpha}.
 \end{eqnarray*}
Take $h= \rho(|x-y|)^{2/\alpha}v_k$, where $v_k$ is a sequence such
that $v_k\rightarrow \infty$ and $v_k^{2-\alpha}k^{-1} \rightarrow
0$. Then one can check that
 \begin{eqnarray*}
\int_{0}^{\infty} D_{l_0(x,y,u)}\phi_k(x-y) \mu_0(du)
 \le
k^{-1}(2-\alpha)^{-1}v_k^{2-\alpha} +(\alpha-1)^{-1}v_k^{1-\alpha},
 \end{eqnarray*}
which tends to zero as $k\rightarrow \infty$. Now since all the
properties in Proposition~\ref{pr3.1} are satisfied we get the pathwise
uniqueness for~(\ref{4.1}). The existence of the solution follows by a
modification of the proof of Theorem~\ref{th2.2}. That gives part (i) of
the theorem. The proof of part (ii) can be given in a similar way. \qed

\bcorollary\label{co4.3} Let $a\ge 0$, $b\ge 0$, $c\ge 0$, $1\le r\le 2$,
$1<\alpha<2$, $q\ge 1$ and $\beta$ be constants. Suppose that $\{B(t)\}$,
$\{L_0(t)\}$ and $\{L_1(t)\}$ are given as the above with $\nu_0(dz) =
z^{-1-\alpha}dz$. If $1/q + 1/\alpha\geq 1$, then for any given $x(0)\in
\mbb{R}_+$ there is a pathwise unique
strong solution to
\begin{eqnarray}\label{4.2}
dx(t) = \sqrt[r]{a |x(t)|} dB(t) + \sign{(x(t-))}\sqrt[q]{c|x(t-)|} dL_0(t) + (\beta x(t)+b)dt + dL_1(t),
 \end{eqnarray}
 and this solution is non-negative.
\ecorollary

\proof One can choose $\rho(z)=\sqrt[r]{z}$ and $p=1/q$ in (4.c) and
hence by Theorem~\ref{th4.2}, there is a pathwise unique strong
solution to (\ref{4.2}) which is non-negative. \qed

In the special case where $r=2$ and $q = \alpha$, the solution of
(\ref{4.2}) is a continuous state branching process with immigration and
the strong existence and uniqueness for (\ref{4.2}) were obtained in
\cite{FL10}.

 \bremark\label{re4.4}
Theorem~\ref{th1.1} follows immediately from Theorem~\ref{th4.2}.
 \eremark

\end{document}